\documentstyle[12pt]{article}
\newcommand{\text}{\mbox}
\begin {document}
\begin{center}
{\large\bf A Remark on ``Two-Sided" Monotonicity Condition: \\
An Application to $L^{p}$ Convergence}\\
\vspace{0.05in}

 S. P. Zhou\footnote{corresponding author} and R. J. Le
\thanks{ Supported in part by Natural Science Foundation of
China under grant number 10471130.}
\end{center}

         \begin{quote}
          \small \bf ABSTRACT.
\rm To verify the universal validity of the ``two-sided" monotonicity condition introduced in [8], we will   apply it to include more classical examples. The present paper selects the $L^{p}$ convergence case for this purpose. Furthermore, Theorem 3 shows that our improvements are not trivial.

          \end{quote}

\begin{center}
\end{center}

\rm  In Fourier analysis, since Fourier coefficients are computable and applicable, people have established many nice results by assuming monotonicity of the coefficients. Generally speaking, it became an important topic how to generalize monotonicity. In many studies  the generalization follows by this way (see, for example, [8] for definitions):
$$\mbox{\rm (coefficients)}\;\;\mbox{\rm nonincreasing}\;\;\Rightarrow\;\;\mbox{\rm quasimonotone}$$
$$\Rightarrow\;\;\mbox{\rm regularly varying quasimonotone}$$
$$\Rightarrow\;\;\mbox{\rm $O$-regularly varying quasimonotone}$$

 \vspace{3mm}

\rm  On the other hand, some mathematicians such as Leindler introduced ``\it rest bounded variation\rm" condition which also generalizes monotonicity: a nonnegative sequence $\{b_{n}\}$ with $b_{n}\to 0$ as $n\to\infty$ is called of ``rest bounded variation" (in symbol, $\{b_{n}\}\in\mbox{\rm $R^{+}_{0}$BV}$) if   $$\sum_{n=m}^{\infty}|b_{n}-b_{n+1}|\leq M(\mbox{\bf b})b_{m} \hspace{.4in}(1)$$
for some constant $M(\mbox{\bf b})$ depending only upon $\mbox{\bf b}$ and $m=1, 2, \cdots$.

\vspace{3mm}

\rm  Since quasimonotonicity  and ``\it rest bounded variation\rm"  are not comparable (cf. [6, Theorem 1]), we suggested the following condition (see [8]) to include both:

\vspace{3mm}

\bf Definition. \it Let \mbox{\bf c}$=\{c_{n}\}_{n=1}^{\infty}$ be a nonnegative sequence tending to zero. If
$$\sum_{n=m}^{2m}|\Delta c_{n}|\leq M(\mbox{\bf c})c_{m}$$
holds for all $m=1, 2, \cdots$, where  $M(\mbox{\bf c})$ indicates a positive constant only depending upon $\mbox{\bf c}$, then we say that the sequence $\mbox{\bf c}$ belongs to class \mbox{\rm GBV}.

\vspace{3mm}

\rm We can verify that either $\{b_{n}\}$ is ($O$-regularly varying) quasimonotone or $\{b_{n}\}\in\mbox{\rm $R^{+}_{0}$BV}$ does imply that $\{b_{n}\}\in \mbox{\rm GBV}$ (Zhou and Le [8, Theorem 3]). The converse is not true (since quasimonotonicity  and ``rest bounded variation" condition (1) are not comparable (cf. [6, Theorem 1])).

 We give more details here. In any sense, monotonicity, quasi-monotonicity, and ``rest bounded variation" condition are all ``one-sided" monotonicity condition, that means, a positive sequence $\mbox{\bf b}=\{b_{n}\}$ under any of these conditions satisfies $b_{n}\leq Cb_{k}$ for $n\geq k$: $b_{n}$ can be majorized by one factor $b_{k}$. But for $\{b_{n}\}\in  \mbox{\rm GBV}$, one can calculate, for $k\leq n\leq 2k$,
$$b_{n}=\sum_{j=n}^{2k}\Delta b_{j}+b_{2k+1}\leq \sum_{j=n}^{2k}|\Delta b_{j}|+b_{2k+1}\leq\sum_{j=k}^{2k}|\Delta b_{j}|+b_{2k+1}\leq M(\mbox{\bf b})b_{k}+b_{2k+1},\;\;\;\;(2)$$
and this can be actually regarded as a ``two-sided" monotonicity: $b_{n}$ is  majorized not only by $b_{k}$, $k\leq n$, but also by $b_{2k+1}$, $2k+1>n$. Therefore, the essential point of \mbox{\rm GBV} condition is to extend monotonicity from ``one-sided" to ``two-sided".

We already reestablished some important results concerning uniform convergence, $L^{1}$-convergence and best approximation rate of certain trigonometric (Fourier) series under GBV condition in [8], [4] and [9] respectively.

To verify the universal validity of this ``two-sided" monotonicity, we need to apply it to include more classical examples. The present paper selects the $L^{p}$ convergence case  for this purpose.

Let $L^{p}$, $1<p<\infty$, be the space of all $p$ power integrable functions of period $2\pi$, $\omega(f, t)_{L^{p}}$ the modulus of continuity in $L^{p}$ norm. Write
$$f(x)=\sum_{n=0}^{\infty}b_{n}\cos nx,\;\;\;\; g(x)=\sum_{n=0}^{\infty}a_{n}\sin nx$$
at any point $x$ where the series converges. Denote $\phi(x)$ to stand for $f$ or $g$ and the $\lambda_{n}$ are its associated Fourier coefficients.

The first group of results is about the necessary and sufficient condition of a sum function beloning to $L^{p}$ spaces.

\vspace{3mm}

\bf Theorem 1.1 ([2] or [3, p.37]). \it Let $1<p<\infty$. If $\lambda_{n}\geq 0$ and $1/p-1<\gamma <1/p$, then a sufficient condition for $\phi(x)x^{-\gamma}\in L^{p} $ is
$$\sum_{n=1}^{\infty}n^{p+p\gamma-2}\left(\sum_{k=n}^{\infty}|\lambda_{k}-\lambda_{k+1}|\right)^{p}<\infty,\hspace{.4in}(3)$$
and a necessary condition is
$$\sum_{n=1}^{\infty} n^{p+p\gamma-2}\left(\sum_{k=n}^{\infty}k^{-1}\lambda_{k}\right)^{p}<\infty.\hspace{.4in}(4)$$

\bf Theorem 1.2 ([6, Theorem 4]). \it If $1<p<\infty$ and $\lambda=\{\lambda_{n}\}\in \mbox{\rm $R^{+}_{0}$BV}$, then $x^{-\gamma}\phi(x)\in L^{p}$, $1/p-1<\gamma <1/p$, if and only if
$$\sum_{n=1}^{\infty} n^{p+p\gamma-2}\lambda_{n}^{p}<\infty.$$

\rm Our result is the following

\vspace{3mm}

\bf Theorem 1.3. \it If $1<p<\infty$ and $\lambda=\{\lambda_{n}\}\in \mbox{\rm GBV}$, then $x^{-\gamma}\phi(x)\in L^{p}$, $1/p-1<\gamma <1/p$, if and only if
$$\sum_{n=1}^{\infty} n^{p+p\gamma-2}\lambda_{n}^{p}<\infty.\hspace{.4in}(5)$$

\rm We make a remark here. Boas in [3] raised Question 6.12: ``What condition is necessary and sufficient for $x^{-\gamma}\phi(x)\in L^{p}$, $1/p-1<\gamma<1/p$, when $\lambda_{n}\geq 0$?"

As we can see from the example given in Theorem 3 (originally given by Leindler [6]), Theorem 1.3 does give a quite gratifying answer up to date under the quite natural GBV condition.

\rm Next, we consider to generalize the $L^{p}$ convergence rate.

\vspace{3mm}

\bf Theorem 2.1 ([1]). \it Let $\{\lambda_{n}\}$ be a decreasing sequence tending to zero such for a fixed $p$, $1<p<\infty$, that
$$\sum_{n=1}^{\infty}n^{p-2}\lambda_{n}^{p}<\infty.\hspace{.4in}(6)$$
Then
$$\omega(f, n^{-1})_{L^{p}}\leq Cn^{-1}\left(\sum_{k=1}^{n-1}k^{2p-2}\lambda_{k}^{p}\right)^{1/p}+C\left(\sum_{k=n}^{\infty}k^{p-2}\lambda^{p}_{k}\right)^{1/p}.\hspace{.4in}(7)$$

\bf Theorem 2.2 ([7, Theorem 1]). \it Let $\{\lambda_{n}\}\in \mbox{\rm $R^{+}_{0}$BV}$ such for a fixed $p$, $1<p<\infty$, that $(6)$ holds,
 then $(7)$ follows.

\vspace{3mm}

\rm Our result can be read as

\vspace{3mm}

\bf Theorem 2.3. \it Let $\{\lambda_{n}\}\in \mbox{\rm GBV}$ such for a fixed $p$, $1<p<\infty$, that $(6)$ holds,
 then $(7)$ follows.

\vspace{3mm}

\rm The following theorem shows that our results do give essential improvements to the previous.

\vspace{3mm}

\bf Theorem 3. \it There exists a nonnegative sequence $\{\lambda_{n}\}\in\mbox{\rm GBV}$ which does not belong to $\mbox{\rm $R^{+}_{0}$BV}$ such that
$$\sum_{n=1}^{\infty} n^{p+p\gamma-2}\lambda_{n}^{p}<\infty$$
for $1<p<\infty$ and $1/p-1<\gamma <1/p$.

\vspace{3mm}
\bf Lemma 1 ([8, Theorem 3]). \it If $\{\lambda_{n}\}$ is an $O$-regularly varying quasimonotone sequence, then  $\{\lambda_{n}\}\in\mbox{\rm GBV}$.

\vspace{3mm}
\bf Proof of Theorem 3.

\vspace{3mm}
\rm We check the example introduced by Leindler [6]: denote $v_{m}=2^{2^{m}}$, and
$$\lambda_{n}=\left\{\begin{array}{ll}
1/(m^{2}v_{m+1})&\mbox{\rm if}\;\;n=v_{m},\\
d_{v_{m}}\Pi_{l=v_{m}}^{n-1}(1+1/l)&\mbox{\rm if}\;\;v_{m}+1\leq n\leq mv_{m},\\
d_{mv_{m}}&\mbox{\rm if}\;\;mv_{m}<n<v_{m+1}.\end{array}\right.$$
Leindler already proved that $\{\lambda_{n}\}$ is a quasi-monotone sequence (then it is in GBV by Lemma 1) but it does not belong to $\mbox{\rm $R^{+}_{0}$BV}$. He pointed out that it is clear that $\lim\limits_{n\to\infty}n\lambda_{n}=0$. Therefore we have
$$\sum_{n=1}^{\infty} n^{p+p\gamma-2}\lambda_{n}^{p}\leq\sum_{n=1}^{\infty} n^{p\gamma-2}(n\lambda_{n})^{p}\leq\sum_{n=1}^{\infty}n^{p\gamma-2}<\infty$$
if considering $\gamma <1/p$. \vspace{3mm}

\rm Throughout the paper, we always use $C$ to indicate an absolute positive constant which may have different values at different occurrences if not specified.

\vspace{3mm}
 \bf Lemma 2 ([5, Theorem 1]). \it If $p\geq 1$ and $\alpha_{n}\geq 0$, then for any sequence $\{\mu_{n}\}$ of positive numbers, it holds that
$$\sum_{n=1}^{\infty}\mu_{n}\left(\sum_{k=1}^{n}\alpha_{k}\right)^{p}\leq p^{p}\sum_{n=1}^{\infty}\mu_{n}^{1-p}\left(\sum_{k=n}^{\infty}\mu_{k}\right)^{p}\alpha_{n}^{p},\hspace{.4in}(8)$$
and
$$\sum_{n=1}^{\infty}\mu_{n}\left(\sum_{k=n}^{\infty}\alpha_{k}\right)^{p}\leq p^{p}\sum_{n=1}^{\infty}\mu_{n}^{1-p}\left(\sum_{k=1}^{n}\mu_{k}\right)^{p}\alpha_{n}^{p}.\hspace{.4in}(9)$$

\bf Lemma 3. \it Let  $\{\lambda_{n}\}\in\mbox{\rm GBV}$, then for $n\geq 1$,
$$\sum_{j=1}^{\infty} \lambda_{2^{j}n}\leq C\sum_{k=n}^{\infty}k^{-1}\lambda_{k}.$$

\bf Proof. \rm From (2), we see that, for $k\leq m\leq 2k$,
$$\lambda_{m}\leq M(\mbox{\bf $\lambda$})\lambda_{k}+\lambda_{2k+1}, $$
thus for $2^{j-1}n\leq k\leq 2^{j}n-1$, $2^{j}n\leq 2k+1\leq 2^{j+1}n-1$, $j=1, 2, \cdots$, we have
$$\lambda_{2^{j}n}\leq M(\mbox{\bf $\lambda$})(\lambda_{k}+\lambda_{2k+1}),$$
or
$$\lambda_{2^{j}n}\leq M(\mbox{\bf $\lambda$})2^{-j+1}n^{-1}\sum_{k=2^{j-1}n}^{2^{j}n-1}( \lambda_{k}+\lambda_{2k+1})$$
$$\leq M(\mbox{\bf $\lambda$})2^{-j+1}n^{-1}\sum_{k=2^{j-1}n}^{2^{j+1}n-1}\lambda_{k}\leq 4M(\mbox{\bf $\lambda$})\sum_{k=2^{j-1}n}^{2^{j+1}n-1}k^{-1}\lambda_{k}.$$
Summing up all the terms from $j=1, 2, \cdots$, we achieve the
required result.

\vspace{3mm}

\bf Lemma 4. \it Let $1<p<\infty$, $\lambda_{n}\geq 0$, and $\sum_{n=1}^{\infty}n^{p-2}\lambda_{n}^{p}<\infty$. Then for $n\geq 1$,
$$n^{1-1/p}\sum_{k=[n/2]+1}^{\infty}k^{-1}\lambda_{k}\leq Cn^{-1}\left(\sum_{k=1}^{n-1}k^{2p-2}\lambda_{k}^{p}\right)^{1/p}+ C\left(\sum_{k=n}^{\infty}k^{p-2}\lambda_{k}^{p}\right)^{1/p}.$$

\bf Proof. \rm  Write
$$n^{p-1}\left(\sum_{j=[n/2]+1}^{\infty}j^{-1}\lambda_{j}\right)^{p}\leq Cn^{-p}\sum_{k=1}^{[n/2]-1}k^{2p-2}\left(\sum_{j=[n/2]+1}^{\infty}j^{-1}\lambda_{j}\right)^{p}$$
$$\leq Cn^{-p}\sum_{k=1}^{[n/2]-1}k^{2p-2}\left(\sum_{j=k}^{\infty}j^{-1}\lambda_{j}\right)^{p}.$$
Put
$$\beta_{k}=\left\{\begin{array}{ll}
k^{2p-2},&k<n,\\
n^{2p}k^{-2},& k\geq n.\end{array}\right.$$
Then
$$ n^{p-1}\left(\sum_{j=[n/2]+1}^{\infty}j^{-1}\lambda_{j}\right)^{p}\leq n^{-p}\sum_{k=1}^{\infty}\beta_{k}\left(\sum_{j=k}^{\infty}j^{-1}\lambda_{j}\right)^{p}.$$
Applying (9), we get
$$\sum_{k=1}^{\infty}\beta_{k}\left(\sum_{j=k}^{\infty}j^{-1}\lambda_{j}\right)^{p}\leq p^{p}\sum_{k=1}^{\infty}\beta_{k}^{1-p}\left(\sum_{j=1}^{k}\beta_{j}\right)^{p}k^{-p}\lambda_{k}^{p}$$
$$= p^{p}\sum_{k=1}^{n-1}k^{(2p-2)(1-p)} \left(\sum_{j=1}^{k}j^{2p-2}\right)^{p}k^{-p}\lambda_{k}^{p}$$
$$+ p^{p}\sum_{k=n}^{\infty}n^{2p(1-p)}k^{-2(1-p)} \left(\sum_{j=1}^{n-1} j^{2p-2}+\sum_{j=n}^{k} n^{2p}j^{-2}\right)^{p}k^{-p}\lambda_{k}^{p}=:J_{1}+J_{2}.$$
It is clear that
$$J_{1}\leq C^{p}\sum_{k=1}^{n-1}k^{2p-2}\lambda_{k}^{p}.$$
At the same time,
$$\sum_{k=n}^{\infty}n^{2p(1-p)}k^{-2(1-p)} \left(\sum_{j=1}^{n-1} j^{2p-2}\right)^{p}k^{-p}\lambda_{k}^{p}\leq Cn^{p}\sum_{k=n}^{\infty}k^{p-2}\lambda_{k}^{p},$$
while
$$\sum_{k=n}^{\infty}n^{2p(1-p)}k^{-2(1-p)} \left(\sum_{j=n}^{k} n^{2p}j^{-2}\right)^{p}k^{-p}\lambda_{k}^{p}\leq Cn^{p}\sum_{k=n}^{\infty}k^{p-2}\lambda_{k}^{p},$$
and therefore,
$$J_{2}\leq C^{p}n^{p}\sum_{k=n}^{\infty}k^{p-2}\lambda_{k}^{p}.$$
Altogether, we have the required inequality.

\vspace{3mm}

\bf Proof of Theorem 1.3.
\vspace{3mm}
Sufficiency. \rm Suppose that (5) holds. Then, since $\{\lambda_{n}\}\in\mbox{\rm GBV}$,
$$\sum_{n=1}^{\infty}n^{p+p\gamma-2}\left(\sum_{k=n}^{\infty}|\lambda_{k}-\lambda_{k+1}|\right)^{p}=\sum_{n=1}^{\infty}n^{p+p\gamma-2}\left(\sum_{j=0}^{\infty}\sum_{k=2^{j}n}^{2^{j+1}n-1}|\lambda_{k}-\lambda_{k+1}|\right)^{p} $$
$$\leq M^{p}({\bf\lambda})\sum_{n=1}^{\infty}n^{p+p\gamma-2}\left(\lambda_{n}^{p}+\left(\sum_{j=1}^{\infty} \lambda_{2^{j}n}\right)^{p}\right)$$
$$\leq M^{p}({\bf\lambda})\left(\sum_{n=1}^{\infty}n^{p+p\gamma-2}\lambda^{p}_{n}+\sum_{n=1}^{\infty}n^{p+p\gamma-2}\left(\sum_{k=n}^{\infty}k^{-1}\lambda_{k}\right)^{p}\right)\;\;\mbox{\rm (by Lemma 3)}$$
$$ \leq M^{p}({\bf\lambda})\left(\sum_{n=1}^{\infty}n^{p+p\gamma-2}\lambda^{p}_{n}+\sum_{n=1}^{\infty}n^{(1-p)(p+p\gamma-2)}\left(\sum_{k=1}^{n}k^{p+p\gamma-2}\right)^{p}n^{-p}\lambda^{p}_{n}\right)\;\;\mbox{\rm (by (9))}$$
$$ \leq M^{p}({\bf\lambda})p^{p}\sum_{n=1}^{\infty}n^{p+p\gamma-2}\lambda^{p}_{n}<\infty.$$
By (3) of Theorem 1.1, it follows that $x^{-\gamma}\phi(x)\in L^{p}$.
\vspace{3mm}

\bf Necessity. \rm  If $x^{-\gamma}\phi(x)\in L^{p}$, then (4) holds. We check that,
$$\lambda_{2n}\leq \sum_{k=2n}^{\infty}|\Delta\lambda_k|=\sum_{j=1}^{\infty}\sum_{k=2^{j}n}^{2^{j+1}n-1}|\Delta\lambda_k|\leq M({\bf\lambda})\sum_{j=1}^{\infty}\lambda_{2^{j}n}$$
$$\leq M({\bf\lambda})\sum_{k=n}^{\infty}k^{-1}\lambda_{k}.\;\;\mbox{\rm (by Lemma 3)}$$
Similarly,
$$\lambda_{2n+1}\leq\lambda_{2n}+\sum_{k=2n}^{\infty}|\Delta\lambda_{k}|\leq M({\bf\lambda})\sum_{k=n}^{\infty}k^{-1}\lambda_{k}.$$
Hence, it yields immediately that
$$\sum_{n=2}^{\infty} n^{p+p\gamma-2}\lambda_{n}^{p}
=\sum_{n=1}^{\infty}(2n) ^{p+p\gamma-2}\lambda_{2n}^{p}+\sum_{n=1}^{\infty}(2n+1) ^{p+p\gamma-2}\lambda_{2n+1}^{p}$$
$$\leq M({\bf\lambda})\left(\sum_{n=1}^{\infty} (2n)^{p+p\gamma-2} \left(\sum_{k=n}^{\infty}k^{-1}\lambda_{k}\right)^{p}+\sum_{n=1}^{\infty} (2n+1)^{p+p\gamma-2} \left(\sum_{k=n}^{\infty}k^{-1}\lambda_{k}\right)^{p}\right)$$
$$\leq M^{p}({\bf\lambda})\sum_{n=1}^{\infty} n^{p+p\gamma-2} \left(\sum_{k=n}^{\infty}k^{-1}\lambda_{k}\right)^{p}<\infty.\;\;\Box$$

\bf Proof of Theorem 2.3.
\vspace{3mm}
\rm The condition (6) is the case $\gamma=0$ in condition (5). Following Leindler's basic technique in the proof of Theorem 2.2 in [8], we have
$$\omega(f, \pi/(2n))_{L^{p}}\leq C\sup_{0<t\leq\pi/(2n)}\left\{\left(\int_{0}^{\pi/n}\left|\sum_{k=1}^{n-1}\lambda_{k}\sin\frac{1}{2}kt\sin k(x\pm t/2)\right|^{p}dx\right)^{1/p}\right.$$
$$+ \left(\int_{0}^{\pi/n}\left|\sum_{k=n}^{\infty}\lambda_{k}(\cos k(x\pm t)-\cos kt)\right|^{p}dx\right)^{1/p}$$
$$+ \left(\int_{\pi/n}^{\pi}\left|\sum_{k=1}^{n}\Delta\lambda_{k}(D_{k}(x\pm t)-D_{k}(x))\right|^{p}dx\right)^{1/p}$$
$$+ \left.\left(\int_{\pi/n}^{\pi}\left|\sum_{k=n+1}^{\infty}\Delta\lambda_{k}(D_{k}(x\pm t)-D_{k}(x))\right|^{p}dx\right)^{1/p}\right\}$$
$$=: C\sup_{0<t\leq\pi/(2n)} (I_{11}+I_{12}+I_{21}+I_{22}),\hspace{.4in}(10)$$
where $D_{k}(x)$ is the Dirichlet kernel of order $k$. Without any change as in [8] one still calculate that
$$I_{11}\leq Cn^{-1}\left(\sum_{k=1}^{n-1}k^{2p-2}\lambda_{k}^{p}\right)^{1/p}.\hspace{.4in}(11)$$
With the same idea as in [8], one reach that
$$I_{12}\leq C\left(\sum_{m=n}^{\infty}\int_{3\pi/(2(m+1))}^{ 3\pi/(2m)} \left|\sum_{k=n}^{\infty}\lambda_{k}\cos kx\right|^{p}dx\right)^{1/p}.$$
By different calculation we now proceed that, for $\{\lambda_{n}\}\in\mbox{\rm GBV}$ and $x\in (3\pi/(2(m+1)), 3\pi/(2m))$, by Abel's transformation,
$$\left|\sum_{k=n}^{\infty}\lambda_{k}\cos kx\right|\leq\sum_{k=n}^{m}\lambda_{k}+C(m+1)\lambda_{m+1}+C(m+1)\sum_{j=0}^{\infty}\sum_{k=2^{j}(m+1)}^{2^{j+1}(m+1)-1}|\Delta\lambda_{k}|$$
$$\leq \sum_{k=n}^{m}\lambda_{k}+C(m+1)\lambda_{m+1}+ M({\bf\lambda})(m+1)\sum_{j=1}^{\infty}\lambda_{2^{j}(m+1)}$$
$$ \leq \sum_{k=n}^{m}\lambda_{k}+C(m+1)\lambda_{m+1}+ M({\bf\lambda})(m+1)\sum_{k=m+1}^{\infty}k^{-1}\lambda_{k},\;\;\mbox{\rm (by Lemma 3)}$$
therefore
$$I_{12}^{p}\leq C^{p}\sum_{m=n}^{\infty}m^{-2}\left(\sum_{k=n}^{m}\lambda_{k}\right)^{p}+C^{p}\sum_{m=n}^{\infty}m^{p-2}\lambda_{m}^{p}+ M^{p}({\bf\lambda})\sum_{m=n}^{\infty}m^{p-2}\left(\sum_{k=m}^{\infty}k^{-1}\lambda_{k}\right)^{p}.$$
By using the inequality (8) of Lemma 2, setting $\mu_{m}=m^{-2}$,
$\alpha_{m}=0$ for $m<n$ and $\alpha_{m}=\lambda_{m}$ otherwise,
we see that
$$\sum_{m=n}^{\infty}m^{-2}\left(\sum_{k=1}^{m}\alpha_{k}\right)^{p}\leq\sum_{m=1}^{\infty}m^{-2}\left(\sum_{k=1}^{m}\alpha_{k}\right)^{p}\leq C^{p}\sum_{m=1}^{\infty}m^{p-2}\alpha_{m}^{p}= C^{p}\sum_{m=n}^{\infty}m^{p-2}\lambda_{m}^{p}.$$
At the same time, by a similar argument, with (9) instead, it yields that  $$\sum_{m=n}^{\infty}m^{p-2}\left(\sum_{k=m}^{\infty}k^{-1}\lambda_{k}\right)^{p}=\sum_{m=1}^{\infty}(n+m-1)^{p-2}\left(\sum_{k=m}^{\infty}(n+k-1)^{-1}\lambda_{n+k-1}\right)^{p}$$
$$\leq C^{p}\sum_{m=1}^{\infty}(n+m-1)^{p-2}\lambda_{n+m-1}^{p}=C^{p}\sum_{m=n}^{\infty}m^{p-2}\lambda_{m}^{p}.$$
Altogether,
$$I_{12}\leq M(\mbox{\bf $\lambda$})\left(\sum_{m=n}^{\infty}m^{p-2}\lambda_{m}^{p}\right)^{1/p}.\hspace{.4in}(12)$$
Following the proof of Theorem 2.1 (see also [7, (4.2)]) we get
$$I_{21}^{p}\leq C^{p}n^{-p}\left(\sum_{m=1}^{n-1}m^{-2}\left(\sum_{k=1}^{m}k^{2}|\Delta\lambda_{k}|\right)^{p}+  \left(\sum_{m=1}^{n-1}m^{p-2}\left(\sum_{k=m+1}^{n}k |\Delta\lambda_{k}|\right)^{p}\right)\right).$$
By the standard technique as we discussed above (by a similar method to Lemma 3), we then have (in view of that $\{\lambda_{n}\}\in\mbox{\rm GBV}$)
$$\sum_{k=1}^{m}k^{2}|\Delta\lambda_{k}|\leq M({\bf\lambda})\sum_{k=1}^{m}k\lambda_{k},$$
and
$$\sum_{k=m+1}^{n}k |\Delta\lambda_{k}|\leq M({\bf\lambda})\left(\sum_{k=m+1}^{n}\lambda_{k}+m\lambda_{m}\right).$$
These are the same as Leindler did in [7], thus the same estimate can be read as
$$I_{21}\leq M({\bf\lambda})n^{-1}\left(\sum_{k=1}^{n-1}k^{2p-2}\lambda_{k}^{p}\right)^{1/p}.\hspace{.4in}(13)$$
Finally we calculate $I_{22}$. Check again that, for $\{\lambda_{n}\}\in\mbox{\rm GBV}$,
$$\sum_{k=n+1}^{\infty}|\Delta\lambda_{k}|=\sum_{j=0}^{\infty}\sum_{k=2^{j}(n+1)}^{2^{j+1}(n+1)-1} |\Delta\lambda_{k}|\leq \lambda_{n+1}+M({\bf\lambda})\sum_{j=1}^{\infty}\lambda_{2^{j}(n+1)}$$
$$\leq \lambda_{n+1}+M({\bf\lambda})\sum_{k=n+1}^{\infty}k^{-1}\lambda_{k}.$$
In a similar way as we did in the proof of Lemma 3, it yields that
$$\lambda_{n+1}\leq C\sum_{k=[n/2]+1}^{2n}k^{-1}\lambda_{k}.$$
Now
$$I_{22}\leq 2\left(\int_{\pi/(2n)}^{\pi+\pi/(2n)}\left(\sum_{k=n+1}^{\infty}|\Delta\lambda_{k}||D_{k}(x)|\right)^{p}dx\right)^{1/p}$$
$$\leq M({\bf\lambda})\sum_{k=[n/2]+1}^{\infty}k^{-1}\lambda_{k}\left(\int_{\pi/(2n)}^{\infty}x^{-p}dx\right)^{1/p}\leq M({\bf\lambda})n^{1-1/p}\sum_{k=[n/2]+1}^{\infty}k^{-1}\lambda_{k},$$
so that\footnote{We remark that this estimate is different from that  Leindler did in [8], where he estimated that $I_{22}\leq Cn^{-1}
\left(\sum_{k=1}^{n-1}k^{2p-2}\lambda_{k}^{p}\right)^{1/p}$.}
$$ I_{22}\leq M({\bf\lambda})n^{-1} \left(\sum_{k=1}^{n-1}k^{2p-2}\lambda_{k}^{p}\right)^{1/p}+ M({\bf\lambda}) \left(\sum_{k=n}^{\infty}k^{p-2}\lambda_{k}^{p}\right)^{1/p}\hspace{.4in}(14)$$
by Lemma 4. With all the estimates (10)-(14), we reach the
required result.

\vspace{3mm}

\begin{center}
{\Large\bf References}
\end{center}
\begin{enumerate}
\rm
\item S. Aljancic, \it On the integral moduli of continuity in $L_{p}$ $(1<p<\infty)$ of Fourier series with monotone coefficients, \rm Proc. Amer. Math. Soc., 17(1966), 287-294.
\item R. A. Askey and S. Wainger, \it Integrability theorems for Fourier series, \rm Duke Math. J., 33(1966), 223-228.
\item R. P. Boas, Jr., \it Integrability Theorems for Trigonometric Transforms, \rm Springer, Berlin-Heidelberg, 1967.
\item R. J. Le and S. P. Zhou,  \it On $L^{1}$ convergence of Fourier series of complex valued functions, \rm Studia Sci. Math. Hungar., tp appear.
\item L. Leindler, \it Generalization of inequalities of Hardy and Littlewood, \rm Acta Sci. Math. (Szeged), 31(1970), 279-285.
\item L. Leindler, \it A new class of numerical sequences and its applications to sine and cosine series, \rm Anal. Math., 28(2002), 279-286.
  \item L. Leindler, \it Relations among Fourier coefficients and sum-functions, \rm Acta Math. Hungar., 104(2004), 171-183.
\item  R. J. Le and S. P. Zhou, \it A new condition for the uniform convergence of certain trigonometric series, \rm Acta Math. Hungar., 108(2005), 161-169.
\item S. P. Zhou and R. J. Le, \it Best approximation rate of certain trigonometric series: important applications of a ``two-sided" monotonicity condition, \rm Acta Math. Sinica, 49(2006), 509-518. (in Chinese)

\end{enumerate}

\begin{flushleft}
\rm S. P. Zhou:\\
Institute of Mathematics\\
Zhejiang  Sci-Tech University\\
Xiasha Economic Development Area\\
Hangzhou, Zhejiang 310018  China\\
Email: szhou@zjip.com
\end{flushleft}

\begin{flushleft}
\rm R. J. Le:\\
 Department of Mathematics\\
Ningbo University\\
Ningbo, Zhejiang 315211 China\\
Email: le\_ye@nbip.net
\end{flushleft}

\begin{center}
\bf Keywords \rm convergence,  quasimonotone,  bounded variation, ``two-sided" monotonicity
\end{center}
\end{document}